\documentclass[12pt,psamsfonts]{amsart}
\usepackage{amsmath}
\usepackage{amsthm}
\usepackage{amssymb}
\usepackage{amscd}
\usepackage{amsfonts}
\usepackage{amsbsy}
\usepackage{color}
\usepackage{epsfig}
\usepackage{graphicx,psfrag}
\usepackage{url}

\newcommand{\R}{{\mathbb R}}

\newtheorem {theorem} {Theorem} 

\newcommand{\al}{\alpha}
\newcommand{\be}{\beta}
\newcommand{\ga}{\gamma}

\begin{document}

\title[Discontinuous piecewise linear differential systems]
{On the number of limit cycles in discontinuous piecewise linear
differential systems with two pieces separated by a straight line}

\author[R.D. Euz\'{e}bio and J. Llibre]
{Rodrigo D. Euz\'{e}bio$^1$ and Jaume Llibre$^2$}

\address{$^1$ Departament de Matem\'atica,
IBILCE, UNESP, Rua Cristovao Colombo 2265, Jardim Nazareth, CEP
15.054-00, Sao Jos\'e de Rio Preto, SP, Brazil}
\email{rodrigo.euzebio@sjrp.unesp.br}

\address{$^2$ Departament de Matem\`{a}tiques,
Universitat Aut\`{o}noma de Barcelona, 08193 Bellaterra, Barcelona,
Catalonia, Spain} \email{jllibre@mat.uab.cat}



\date{}
\dedicatory{}

\maketitle

\begin{abstract}
In this paper we study the maximum number $N$ of limit cycles that
can exhibit a planar piecewise linear differential system formed by
two pieces separated by a straight line. More precisely, we prove
that this maximum number satisfies $2\le N \le 3$ if one of the two
linear differential systems has its equilibrium point on the
straight line of discontinuity.
\end{abstract}

\section{Introduction and statement of the main result}\label{s0}

The study of piecewise linear differential systems goes back to
Andronov, Vitt and Khaikin \cite{AVK} and still continues to receive
attention by researchers. These last years a renewed interest has
appeared in the mathematical community working in differential
equations for understanding the dynamical richness of the piecewise
linear differential systems, because these systems are widely used
to model processes appearing in electronics, mechanics, economy,
..., see for instance the books of di Bernardo, Budd, Champneys and
Kowalczyk \cite{dBBCKN}, and Simpson \cite{Si}, and the survey of
Makarenkov and Lamb \cite{ML}, and the hundreds of references quoted
in these last three works.

\smallskip

We recall that a {\it limit cycle} is a periodic orbit of a
differential system which is isolated in the set of all periodic
orbits of the system.

\smallskip

The simplest possible continuous but nonsmooth piecewise linear
differential systems are the ones having only two pieces separated
by a straight line. In 1990 Lum and Chua \cite{LC} conjectured that
a continuous piecewise linear vector field in the plane with two
pieces has at most one limit cycle. We note that even in this
apparent simple case, only after a difficult analysis it was
possible to prove the existence of at most one limit cycle, thus in
1998 this conjecture was proved by Freire, Ponce, Rodrigo and Torres
\cite{FPRT}. There are two reasons that difficult the analysis of
these differential systems. First, even one can easily integrate the
solutions of every linear differential system, the time that an
orbit expends in each half--plane governed by each linear
differential system is in general unknown, consequently the matching
of the corresponding solutions is a difficult problem. Second, the
number of parameters to consider in order to be sure that we take
into account all possible cases is in general not small. Of course,
these difficulties increase when we work with discontinuous
piecewise linear differential systems. Recently, a new an easier
proof that at most one limit cycle exists for the continuous
piecewise linear differential systems with two pieces separated by a
straight line has been done by Llibre, Ordo\~{n}ez and Ponce in
\cite{LOP}.

\smallskip

The objective of this paper is to study the problem of Lum and Chua
but now for the class of discontinuous piecewise linear differential
systems in the plane with two pieces separated by a straight line.
In some sense this problem can be seen as an extension of the 16th
Hilbert's problem to the discontinuous piecewise linear differential
systems in the plane with two pieces separated by a straight line.
We recall that the 16th Hilbert's problem essentially ask for the
maximum number of limit cycles that a polynomial differential system
in the plane can have in function of the degree of the system. For
the moment this problem remains open, for more details on the 16th
Hilbert's problem see for instance \cite{Hi, Il, Li}.

\smallskip

Several authors tried to determine the maximum number of nested
limit cycles surrounding a unique equilibrium point for the class of
all discontinuous piecewise linear differential systems with two
pieces separated by a straight line. Thus in the paper of Han and
Zhang \cite{HZ} some results about the existence of two limit cycles
appeared, so that the authors conjectured that the maximum number of
limit cycles for this class of piecewise linear differential systems
is exactly two. But Huan and Yang in \cite{HY} provided numerical
evidence about the existence of three nested limit cycles
surrounding a unique equilibrium. Llibre and Ponce in \cite{LP}
inspired in the numerical example of \cite{HY} proved that there are
discontinuous piecewise linear differential systems with two pieces
separated by a straight line having three limit cycles. Later on
other authors obtained also three limit cycles for those
differential systems following different ways, see the papers of
Braga and Mello \cite{BdM}, of Buzzi, Pessoa and Torregrosa
\cite{BPT}, and of Freire, Ponce and Torres \cite{FPT14}.

\smallskip

The linear differential systems that we consider in every
half--plane extended to the full plane is either a focus (F) (we
include in this class of foci the centers), or a node (N), or a
saddle (S). We recall that there are three classes of linear nodes:
nodes with different eigenvalues, nodes with equal eigenvalues whose
linear part does not diagonalize, and nodes with equal eigenvalues
whose linear part diagonalize, called {\it star nodes}. Clearly if a
piecewise linear differential system with two pieces separated by a
straight line has a star node, this prevents the existence of
periodic orbits.

\smallskip

An equilibrium point $p$ of a linear differential system defined in
a half--plane having in the full plane a node, a focus or a saddle
is {\it real} when $p$ belongs to the closure of the half--plane
where the system is defined the mentioned linear differential
system, and $p$ is called {\it virtual} otherwise.

\smallskip

We distinguish six classes or types of planar discontinuous
piecewise linear differential systems: FF, FN, FS, NN, NS and SS.
Inside these classes and {\it in this paper we only consider limit
cycles surrounding a unique equilibrium point or a unique sliding
segment}. So we do not consider non--sliding limit cycles. Now we
recall the definitions of sliding segment and non--sliding limit
cycle, for more details on these definitions see for instance
\cite{GP} and \cite{Te}.

\smallskip

Let $Z= (X,Y)$ be a discontinuous piecewise linear differential
vector field with two pieces separated by a straight line $\Sigma$,
in one piece we have the linear vector field $X$ and in the other
the linear vector field $Y$. Following Filippov \cite{Fi} we
distinguish three open regions in the discontinuity straight line
$\Sigma$.
\begin{itemize}
\item[1)] The {\it sliding region} $\Sigma^{sl}$ where the vectors
$X(p)$ and $Y(p)$ with $p\in \Sigma$ point inward $\Sigma$.

\item[2)] The {\it escaping region} $\Sigma^e$ where the vectors
$X(p)$ and $Y(p)$ with $p\in \Sigma$ point outward $\Sigma$.

\item[3)] The {\it sewing region} $\Sigma^s$ where the vectors
$X(p)$ and $Y(p)$ with $p\in \Sigma$ point to the same direction and
are transverse to $\Sigma$.
\end{itemize}
Any segment contained in $\Sigma^e\cup \Sigma^{sl}$ is called a {\it
sliding segment}. Any limit cycle $\ga$ of $Z$ such that $\ga \cap
(\Sigma^e \cup  \Sigma^{sl})= \emptyset$ is called a {\it
non--sliding limit cycle}.

\smallskip

Limit cycles of discontinuous piecewise linear differential systems
with two pieces separated by a straight line have been studied by
many authors, see for instance the articles \cite{ALMT, FPT12,
FPT14, HZ, HY, HY1, HY2, LP, LT, SZL}. Summarizing the results of
these articles we have that the maximum number of known limit cycles
that one of these systems can exhibit is given in Table \ref{T1}. In
that table the symbol ``--'' indicates that those cases appear
repeated in the table, because for instance the case NF is the same
as the case FN, which appears with a $3$ in the table. For more
details on Table \ref{T1} see mainly the references \cite{LT} of
Llibre, Teixeira and Torregrosa and \cite{FPT14} of Freire, Ponce
and Torres.

\smallskip

\begin{table}[h]
\begin{center}
\begin{tabular}{ccc}
\begin{tabular}{|c||c|c|c|}
\hline
 & F & N & S \\
\hline \hline
F & 3 & 3 & 3 \\
N & -- & 2 & 2\\
S & -- & -- & 2\\
\hline
\end{tabular}
\end{tabular}\\[0.5cm]
\end{center}
\caption{Lower bounds for the maximum number of limit cycles of
discontinuous piecewise linear differential systems with two pieces
separated by a straight line.}\label{T1}
\end{table}

Our main result is the following.

\begin{theorem}\label{t1}
If one of the linear differential systems has its equilibrium point
on the straight line of separation, then the maximum number $N$ of
limit cycles of discontinuous piecewise linear differential systems
with two pieces separated by a straight line satisfies $2\le N \le
3$.
\end{theorem}

Theorem \ref{t1} is proved in section \ref{s2}. As far as we know it
is the first time that an upper bound of $3$ limit cycles is given
for a class of discontinuous piecewise linear differential systems
with two pieces separated by a straight line.

\smallskip

Under the assumptions of Theorem \ref{t1} we must remark that if the
linear differential system whose equilibrium point is on the
straight line of discontinuity is a saddle or a node, then due to
the existence of the invariant straight lines of the saddle and of
the node such discontinuous piecewise linear differential systems
cannot have periodic solutions, and consequently limit cycles. So we
only must prove Theorem \ref{t1} when the linear differential system
whose equilibrium point is on the straight line of discontinuity is
a focus or a center.

\smallskip

There are two main tools in the proof of Theorem \ref{t1}. The first
are the canonical forms of all the possible configurations of the
discontinuous piecewise linear differential systems in the plane
with two pieces separated by a straight line, see section \ref{s01}.
These canonical forms only depend of five parameters and are due to
Freire, Ponce and Torres, see \cite{FPT14}. The second tool are the
extended complete Chebyshev systems, see for more details section
\ref{s1}.

\section{Canonical forms}\label{s01}

We assume without loss of generality that the two pieces in the
plane where are defined the discontinuous piecewise linear
differential systems are the left and the right half--planes
\[
S^-= \{(x,y)\in \R^2 : x\le 0 \}, \qquad S^+= \{(x,y)\in \R^2 : x\ge
0 \}.
\]
Consequently $x=0$ is the straight line of separation between the
two linear differential systems
\begin{equation}\label{--}
\left(
\begin{array}{c}
\dot x\\
\dot y
\end{array}
\right)= \left(
\begin{array}{cc}
a_{11}^- & a_{12}^-\\
a_{21}^- & a_{22}^-
\end{array}
\right)\left(
\begin{array}{c}
x\\
y
\end{array}
\right)+ \left(
\begin{array}{c}
b_1^-\\
b_2^-
\end{array}
\right),
\end{equation}
defined for the $(x,y)\in S^-$, and
\begin{equation}\label{++}
\left(
\begin{array}{c}
\dot x\\
\dot y
\end{array}
\right)= \left(
\begin{array}{cc}
a_{11}^+ & a_{12}^+\\
a_{21}^+ & a_{22}^+
\end{array}
\right)\left(
\begin{array}{c}
x\\
y
\end{array}
\right)+ \left(
\begin{array}{c}
b_1^+\\
b_2^+
\end{array}
\right),
\end{equation}
defined for the $(x,y)\in S^+$. Note that both systems together
depend on twelve parameters.

\smallskip

Now we consider the discontinuous piecewise linear differential
systems
\begin{equation}\label{---}
\left(
\begin{array}{c}
\dot x\\
\dot y
\end{array}
\right)= \left(
\begin{array}{cc}
2\ell & -1\\
\ell^2-\al^2 & 0
\end{array}
\right)\left(
\begin{array}{c}
x\\
y
\end{array}
\right)+ \left(
\begin{array}{c}
0\\
a
\end{array}
\right),
\end{equation}
defined for the $(x,y)\in S^-$, and
\begin{equation}\label{+++}
\left(
\begin{array}{c}
\dot x\\
\dot y
\end{array}
\right)= \left(
\begin{array}{cc}
2r & -1\\
r^2-\be^2 & 0
\end{array}
\right)\left(
\begin{array}{c}
x\\
y
\end{array}
\right)+ \left(
\begin{array}{c}
b\\
c
\end{array}
\right),
\end{equation}
defined for the $(x,y)\in S^+$, where $\al,\be \in \{i,0,1\}$. Of
course $i= \sqrt{-1}$. Note that both systems together depend on
five parameters. We remark that if $\al=i$ then the equilibrium
point of system \eqref{---} has eigenvalues $\ell \pm i$, so it is a
focus if $\ell\neq 0$, and a center if $\ell=0$. If $\al=0$ then
system \eqref{---} is a node with eigenvalue $\ell\neq 0$ of
multiplicity $2$ whose linear part does not diagonalize. If $\al=1$
then system \eqref{---} is a saddle with eigenvalues $\ell-1$ and
$\ell+1$ when $|\ell|<1$, and a node with eigenvalues $\ell-1$ and
$\ell+1$ whose linear part diagonalize when $|\ell|>1$.

\smallskip

Let $U$ be an open subset of $\R^2$. We say that the homeomorphism
$h$ between $U$ and its image by $h$ is a {\it topological
equivalence} between the discontinuous piecewise linear differential
system \eqref{--}$+$\eqref{++} and the discontinuous piecewise
linear differential system \eqref{---}$+$\eqref{+++} if $h$ applies
orbits of system \eqref{--}$+$\eqref{++} contained in $U$ into
orbits of system \eqref{---}$+$\eqref{+++} contained in $h(U)$.

\smallskip

{From} Propositions 1 and 2 of \cite{FPT14} it follows that there
exists a topological equivalence between the phase portrait of the
discontinuous piecewise linear differential system
\eqref{--}$+$\eqref{++} and the phase portrait of the discontinuous
piecewise linear differential system \eqref{---}$+$\eqref{+++}
restricted to the orbits that do not have points in common with the
sliding set of these systems. Therefore, since we are interested in
studying the limit cycles of the system \eqref{--}$+$\eqref{++}
which do not intersect its sliding set, it will be sufficient to
study the limit cycles of the system \eqref{---}$+$\eqref{+++}.

\section{Extended Complete Chebyshev systems}\label{s1}

The set of functions $\{f_0,f_1,...,f_n\}$ defined on the interval
$I$ form an \emph{Extended Chebyshev} system on $I$, if and only if
any nontrivial linear combination of these functions has at most $n$
zeros counting their multiplicities and this number is reached.

\smallskip

The set of functions $\{f_0,f_1,...,f_n\}$ is an \emph{Extended
Complete Chebyshev} system or simply an ECT--\emph{system} on $I$ if
and only if for $k=0,1,\ldots,n$, the subset of functions
$\{f_0,f_1,...,f_k\}$ form an Extended Chebyshev system.

\smallskip

For proving that the set of functions $\{f_0,f_1,...,f_n\}$ is an
ECT--system on $I$ it is sufficient and necessary to show that the
Wronskians
\[
W(f_0,...,f_k)(s)=\left|\begin{array}{cccc}
                      f_0(s) &f_1(s) & \cdots & f_k(s)\\
                      f'_0(s) &f'_1(s) & \cdots & f'_k(s) \\
                      \vdots &\vdots & \ddots & \vdots \\
                      f_0^{(k)}(s) &f_1^{(k)}(s) & \cdots & f_k^{(k)}(s)
                    \end{array}\right|\neq 0,
\]
on $I$ for $k=0,1,\ldots,n$. For more details on ECT--system see the
book \cite{JK}.

\section{Proof of Theorem \ref{t1}}\label{s2}

We separate the proof of Theorem \ref{t1} in three cases.

\bigskip

\noindent{\bf Case 1}: The discontinuous piecewise linear
differential system \eqref{---}$+$\eqref{+++} has one real focus on
the discontinuity straight line and another focus, real or virtual.
Then in system \eqref{---}$+$\eqref{+++} we must take $\al=\be=i$.
We note that in this case the sliding segment is $\{ (0,y) : 0<y<b
\}$.

\smallskip

The solution of system \eqref{---} starting at the point
$(x,y)=(0,y_0)$ is
\begin{equation}\label{e1}
\begin{array}{l}
x(t)= \dfrac{e^{\ell t} (a \cos t-(y_0+\ell (a+\ell y_0))
\sin t) -a}{\ell^2+1}, \vspace{0.2cm}\\
y(t)= \dfrac{e^{\ell t} \left(\left(y_0 \ell^2+2 a \ell+y_0\right)
\cos t-\left(a \left(\ell^2-1\right)+\ell \left(\ell^2+1\right)
y_0\right) \sin t\right)-2 a \ell}{\ell^2+1}.
\end{array}
\end{equation}

The solution of system \eqref{+++} starting at the point
$(x,y)=(0,y_0)$ is
\begin{equation}\label{e2}
\begin{array}{rl}
x(t)=& \dfrac{e^{r t} \left(c \cos t-\left(b
\left(r^2+1\right)+y_0+r
(c+r y_0)\right) \sin t\right)-c}{r^2+1}, \vspace{0.2cm}\\
y(t)=& -\dfrac{1}{r^2+1} \Big( b r^2+2 c r+b+e^{r t} \big(\left(c
(r^2-1)+r (r^2+1) (b+y_0)\right) \sin t \\
& -\left(b r^2+y_0 r^2+2 c r+b+y_0\right) \cos t\big)\Big).
\end{array}
\end{equation}

Let $t_-$ be the finite positive time that an orbit of system
\eqref{---} expends inside $S^-$ starting at the point $(0,y_0)$ and
entering in $S^-$ in forward time, and let $-t_-$ be the finite
positive time that an orbit of system \eqref{---} expends inside
$S^-$ starting at the point $(0,y_0)$ and entering in $S^-$ in
backward time. If such orbits do not exist for system \eqref{---},
then system \eqref{---}$+$\eqref{+++} cannot have periodic solutions
and we are done. So we assume that there are orbits for which the
times $t_-$ or $-t_-$ are well defined.

\smallskip

In a similar way let $t_+$ be the finite positive time that an orbit
of system \eqref{+++} expends inside $S^+$ starting at the point
$(0,y_0)$ and entering in $S^+$ in forward time, and let $-t_+$ be
the finite positive time that an orbit of system \eqref{+++} expends
inside $S^+$ starting at the point $(0,y_0)$ and entering in $S^+$
in backward time. Again we assume that there are orbits for which
the times $t_+$ or $-t_+$ are well defined, otherwise the system
\eqref{---}$+$\eqref{+++} cannot have periodic solutions.

\smallskip

Note that if we have a periodic solution of system
\eqref{---}$+$\eqref{+++} for such an orbit the times $t_-$ and
$t_+$ satisfy that $t_- t_+<0$.

\smallskip

Assume that system \eqref{---}$+$\eqref{+++} has a periodic solution
and let $t_-$ and $t_+$ be the times associated to the two pieces of
this periodic solution. Then one of the following sets of three
equations must be satisfied for a such periodic solution, either
\begin{equation}\label{e2a}
x(t_+)= 0, \quad x(-t_-)= 0, \quad y(-t_-)-y(t_+)= 0,
\end{equation}
or
\begin{equation}\label{e2b}
x(-t_+)= 0, \quad x(t_-)= 0, \quad y(t_-)-y(-t_+)= 0.
\end{equation}
Now we shall assume that equations \eqref{e2a} hold. The proof of
Theorem \ref{t1}, in Case 1 when $c=0$ and if equations \eqref{e2b}
hold, is completely analogous to the proof when equations
\eqref{e2a} hold when $c=0$.

\smallskip

Using \eqref{e1} and \eqref{e2} the three equations \eqref{e2a}
become
\begin{equation*}\label{e3}
\begin{array}{rl}
e_1=& e^{r t_+} \left(c \cos t_+-\left(b \left(r^2+1\right)+y_0+r
(c+r y_0)\right) \sin t_+\right)-c, \vspace{0.2cm}\\
e_2=& e^{\ell t_-} (a \cos t_- +(y_0+\ell(a+\ell y_0)) \sin
t_-)-a, \vspace{0.2cm}\\
e_3=& -(1+r^2) \Big(e^{-\ell t_-} \big((2 a \ell +y_0+ \ell^2 y_0)
\cos t_- \vspace{0.2cm}\\
& +\left(a (\ell^2-1)+\ell (\ell^2+1) y_0\right) \sin t_- -2 a
\ell\big)\Big)
\vspace{0.2cm}\\
&+(\ell^2+1) \Big(e^{r t_+} \big(\left(b+ b r^2
+2 c r+  y_0+ r^2 y_0\right) \cos t_+ -b r^2\vspace{0.2cm}\\
& -\left(c \left(r^2-1\right)+r \left(r^2+1\right) (b+y_0)\right)
\sin t_+ -2 c r-b\big)\Big).
\end{array}
\end{equation*}

By assumption one of the two foci must be on the discontinuity line
$x=0$. We assume that the focus of system \eqref{+++} is on $x=0$,
i.e. we take $c=0$. We must also study the case when the focus of
system \eqref{---} is on $x=0$ (i.e. $a=0$), because in the
expression of system \eqref{---}$+$\eqref{+++} both foci do not play
exactly the same role due to the parameter $b$. Now we will do the
study when the focus of system \eqref{+++} is on $x=0$, i.e. $c=0$.
Later on we shall study the case $a=0$.

\smallskip

Taking $c=0$ equation $e_1=0$ becomes
\[
-e^{r t_+}(1+ r^2) (y_0+b ) \sin t_+=0.
\]
If $y_0 +b=0$, then at most there is one periodic solution, because
from this equation $y_0=-b$, and we are done. So we suppose that
$\sin t_+=0$, i.e. $t_+= \pi$. Now solving equation $e_2=0$ with
respect to the variable $y_0$ we get
\[
y_0= a \, \frac{e^{-\ell t_-}-\cos t_- - \ell \sin t_- }{(\ell^2+1)
\sin t_- }.
\]
Therefore the equation $e_3=0$ writes
\begin{equation}\label{e3a}
\begin{array}{rl}
-\left(e^{\pi  r}+1\right) (b \ell^2-a \ell+b)e^{\ell t_-}  + a
\left(e^{\pi  r}-1\right) e^{\ell t_-}  \cot t_- & \vspace{0.2cm}\\
+ a \csc t_- - a e^{\pi r}e^{2 \ell t_-} \csc t_- &=0,
\end{array}
\end{equation}
or equivalently
\begin{equation}\label{e40}
\begin{array}{rl}
-\left(e^{\pi  r}+1\right) (b \ell^2-a \ell+b) f_0(t_-)  + a
\left(e^{\pi  r}-1\right)  f_1(t_-)- &\vspace{0.2cm}\\
+ a f_2(t_-) - a e^{\pi r}  f_3(t_-) &=0,
\end{array}
\end{equation}
where
\begin{equation}\label{e4a}
\begin{array}{l}
f_0(t_-)= e^{\ell t_-}, \vspace{0.2cm}\\
f_1(t_-)=  e^{\ell t_-}  \cot t_-, \vspace{0.2cm}\\
f_2(t_-)= \csc t_-, \vspace{0.2cm}\\
f_3(t_-)= e^{2 \ell t_-} \csc t_-.
\end{array}
\end{equation}

Now we claim that the set of functions $\{f_0,f_1,f_2,f_3\}$ is an
extended complete Chebyshev system, and consequently the system
\eqref{e2a} can have at most $3$ solutions. {From} section \ref{s1}
this number of solutions is reached if the coefficients of the
functions $f_0,f_1,f_2,f_3$ in \eqref{e40}, as functions of the
parameters of the discontinuous piecewise differential system
\eqref{---}$+$\eqref{+++}, are functionally independent, which is
not the case because the coefficients $a(e^{\pi  r}-1)$, $a$ and $a
e^{\pi r}$ are not functionally independent. But since the three
coefficients of $f_0,f_1,f_2$ are functionally independent the
maximum number of limit cycles of system \eqref{---}$+$\eqref{+++}
in Case 1 will be at least $2$.

\smallskip

Once the claim be proved, the proof of Theorem \ref{t1} in Case 1
with $c=0$ will be completed. Now we prove the claim. For proving
that the set of functions $\{f_0,f_1,f_2,f_3\}$ is an ECT--system it
is sufficient and necessary to show that the Wronskians
$W(f_0,...,f_k)(t_-)$ are not zero for $k=0,1,2,3$. Indeed, we have
\[
\begin{array}{rl}
W(f_0)(t_-)=& e^{\ell t_-} \neq 0, \vspace{0.2cm}\\
W(f_0,f_1)(t_-)=& -e^{2 \ell t_-} \csc ^2 t_-\neq 0, \vspace{0.2cm}\\
W(f_0,f_1,f_2)(t_-)=& -e^{2 \ell t_-} \left(\ell^2+1\right) \csc^3 t_-
\neq 0, \vspace{0.2cm}\\
W(f_0,f_1,f_2,f_3)(t_-)=& -2 e^{4 \ell t_-} \ell (1 + \ell^2)^2
\csc^4 t_-.
\end{array}
\]
Therefore, if $\ell\neq 0$ (i.e. if the equilibrium point of system
\eqref{---} is a focus), then the set of functions
$\{f_0,f_1,f_2,f_3\}$ is an extended complete Chebyshev system.

\smallskip

Assume now that $\ell=0$, i.e. the equilibrium point of system
\eqref{---} is a center. Then equation \eqref{e3a} becomes
\[
-b \left(1+e^{\pi  r}\right)+a \left(e^{\pi  r}-1\right) (\cot t_-
-\csc t_-)=0,
\]
or equivalently
\[
-b \left(1+e^{\pi  r}\right)f_0(t_-)+a \left(e^{\pi  r}-1\right)
f_1(t_-)=0,
\]
where
\[
\begin{array}{l}
f_0(t_-)=1, \vspace{0.2cm}\\
f_1(t_-)= \cot t_- -\csc t_-.
\end{array}
\]
The set of functions $\{f_0,f_1\}$ is an ECT--system because the
Wronskians $W(f_0)(t_-)$ and  $W(f_0,f_1)(t_-)$ are not zero,
because
\[
\begin{array}{rl}
W(f_0)(t_-)=&  1\neq 0, \vspace{0.2cm}\\
W(f_0,f_1)(t_-)=& (\cot t_- -\csc t_-) \csc t_-\le -\dfrac12.
\end{array}
\]
In short this completes the proof of Theorem \ref{t1} in the Case 1
when $c=0$.

\smallskip

Assume $a=0$. Then the equation $e_2=0$ reduces to
\[
e^{-\ell t_-} \left(\ell^2+1\right) y_0 \sin (t_-)=0.
\]
If $y_0=0$, then at most there is one periodic solution, and we are
done, but we note that this periodic solution would be a
non--sliding limit cycle and consequently we must not take it into
account. So we suppose that $\sin t_-=0$, i.e. $t_-= \pi$. Now
solving equation $e_1=0$ with respect to the variable $y_0$ we get
\[
y_0= -\frac{b r^2+c r+b-c \cot t_+ +c e^{-r t_+} \csc t_+}{r^2+1}.
\]
Therefore the equation $e_3=0$ writes
\[
\begin{array}{rl}
- \left(1+e^{\ell \pi }\right)(b r^2+c r+b)e^{r t_+} +c (1-e^{\ell
\pi}) e^{r t_+} \cot t_+& \vspace{0.2cm}\\
-c \csc t_+ + c e^{\pi \ell} e^{2 r t_+} \csc t_+&=0,
\end{array}
\]
or equivalently
\begin{equation*}\label{e4}
\begin{array}{rl}
- \left(1+e^{\ell \pi }\right)(b r^2+c r+b)f_0(t_+) +c (1-e^{\ell
\pi}) f_1(t_+)& \vspace{0.2cm}\\-c f_2(t_+) + c e^{\pi \ell}
f_3(t_+)&=0,
\end{array}
\end{equation*}
where
\[
\begin{array}{l}
f_0(t_-)= e^{r t_+}, \vspace{0.2cm}\\
f_1(t_-)= e^{r t_+} \cot t_+, \vspace{0.2cm}\\
f_2(t_-)= \csc t_+, \vspace{0.2cm}\\
f_3(t_-)= e^{2 r t_+} \csc t_+.
\end{array}
\]
These functions $f_0, f_1, f_2, f_3$ coincide with the functions
\eqref{e4a} if we change $r$ by $\ell$, and $t_+$ by $t_-$. So the
rest of the proof of Theorem \ref{t1} in Case 1 with $a=0$ follows
as in the Case 1 with $c=0$. In summary we have proved Theorem
\ref{t1} in Case 1.

\bigskip

\noindent{\bf Case 2}: The discontinuous piecewise linear
differential system \eqref{---}$+$\eqref{+++} has one real focus on
the discontinuity straight line and one real or virtual node outside
the discontinuity straight line with the eigenvalue $\ell\neq 0$ of
multiplicity two and whose linear part does not diagonalize. Then in
system \eqref{---}$+$\eqref{+++} we must take $\al=0$ and $\be= i$.
The case $\al=i$ and $\be= 0$ follows in a similar way.

\smallskip

The solution of system \eqref{---} starting at the point
$(x,y)=(0,y_0)$ is
\begin{equation}\label{e10}
\begin{array}{l}
x(t)= \dfrac{e^{\ell t} (a- \ell^2 y_0 t-a \ell t)-a}{\ell^2}, \vspace{0.2cm}\\
y(t)= -\dfrac{2 a+e^{\ell t} (a (\ell t-2)+\ell (\ell t-1)
y_0)}{\ell}.
\end{array}
\end{equation}

The solution of system \eqref{+++} starting at the point
$(x,y)=(0,y_0)$ is given in \eqref{e2}.

\smallskip

Let $t_-$ and  $t_+$ be the finite positive times defined in a
similar way to the Case 1. Again if we have a periodic solution of
system \eqref{---}$+$\eqref{+++} for such an orbit the times $t_-$
and $t_+$ satisfy that $t_- t_+<0$. Suppose that system
\eqref{---}$+$\eqref{+++} has a periodic solution and let $t_-$ and
$t_+$ the times associated to the two pieces of this periodic
solution. Then one of the following sets of three equations
\eqref{e2a}, or \eqref{e2b} must be satisfied.

\smallskip

Now we shall assume that equations \eqref{e2a} hold. Again the proof
of Theorem \ref{t1} in Case 2 if equations \eqref{e2b} hold is
completely analogous to the proof when equations \eqref{e2a} hold.

\smallskip

Using \eqref{e2}, \eqref{e10}, taking into account that $\ell\neq 0$
and that the focus must be on the discontinuity straight line (i.e.
$c=0$) the three equations \eqref{e2a} become
\begin{equation*}\label{e3}
\begin{array}{rl}
e_1=& -e^{r t_+} (b+y_0) \sin t_+=0, \vspace{0.2cm}\\
e_2=& \ell^2 y_0 t_- +a \ell
t_- -a e^{\ell t_-}+a=0, \vspace{0.2cm}\\
e_3=& \ell y_0(1+ \ell t_-) +b \ell e^{\ell t_-} +a\ell t_-  +2
a(1-e^{\ell t_-}) \vspace{0.2cm}\\
& -(b+y_0)\ell e^{\ell t_- +r t_+} (\cos t_+ - r \sin t_+) =0.
\end{array}
\end{equation*}

{From} equation $e_1=0$ if $y_0 +b=0$, then at most there is one
periodic solution, and we are done. So we suppose that $\sin t_+=0$,
i.e. $t_+= \pi$. Now solving equation $e_2=0$ with respect to the
variable $y_0$ we get
\[
y_0= \frac{a \left(e^{\ell t_-}-\ell t_--1\right)}{\ell^2 t_-}.
\]
Therefore the equation $e_3=0$ writes
\begin{equation*}\label{e11}
a \left(1-e^{\pi  r}\right)+ a e^{\pi  r} e^{\ell t_-} -a e^{-\ell
t_-} -\left(1+e^{\pi  r}\right) \ell (a-b \ell) t_-=0,
\end{equation*}
or equivalently
\begin{equation*}\label{e12}
a \left(1-e^{\pi  r}\right)f_0(t_-)+ a e^{\pi  r} f_1(t_-) -a
f_2(t_-) -\left(1+e^{\pi  r}\right) \ell (a-b \ell)f_3(t_-)=0,
\end{equation*}
where
\begin{equation*}\label{e13}
\begin{array}{l}
f_0(t_-)= 1, \vspace{0.2cm}\\
f_1(t_-)= e^{\ell t_-}, \vspace{0.2cm}\\
f_2(t_-)= e^{-\ell t_-}, \vspace{0.2cm}\\
f_3(t_-)= t_-.
\end{array}
\end{equation*}

The set of functions $\{f_0,f_1,f_2,f_3\}$ is an ECT--system because
the Wronskians $W(f_0,...,f_k)(t_-)$ are not zero for $k=0,1,2,3$.
Indeed, we have
\[
\begin{array}{rl}
W(f_0)(t_-)=& 1 \neq 0, \vspace{0.2cm}\\
W(f_0,f_1)(t_-)=& \ell e^{\ell t_-}\neq 0, \vspace{0.2cm}\\
W(f_0,f_1,f_2)(t_-)=& 2\ell^3 \neq 0, \vspace{0.2cm}\\
W(f_0,f_1,f_2,f_3)(t_-)=& -2\ell^5.
\end{array}
\]
Therefore the system \eqref{e2a} can have at most $3$ solutions.
Since the three coefficients of $f_0,f_1,f_2$ are functionally
dependent the system \eqref{---}$+$\eqref{+++} perhaps do not reach
the three solutions. But since the coefficients of $f_1,f_2,f_3$ are
functionally independent the system \eqref{---}$+$\eqref{+++} can
have two solutions. This completes the proof of Theorem \ref{t1} in
Case 2.

\bigskip

\noindent{\bf Case 3}: The discontinuous piecewise linear
differential system \eqref{---}$+$\eqref{+++} has one real focus on
the discontinuity straight line and one real or virtual node or
saddle outside the discontinuity straight line. This node has two
different eigenvalues. Then in system \eqref{---}$+$\eqref{+++} we
must take $\al=1$ and $\be= i$. The case $\al=i$ and $\be= 1$
follows in a similar way.

\smallskip

We recall that if $|\ell|>1$ then system \eqref{---} has a real or
virtual node with two different eigenvalues, while if $|\ell|<1$ the
system has a real or virtual saddle. Both cases are studied
simultaneously.

\smallskip

The solution of system \eqref{---} starting at the point
$(x,y)=(0,y_0)$ is
\begin{equation}\label{e20}
\begin{array}{rl}
x(t)=& -\dfrac{e^{-t}}{2(\ell^2-1)} \Big(a \left(e^{(\ell+2) t}
(\ell-1)+2 e^t-e^{\ell t}(\ell+1)\right)\\
& \qquad \qquad \qquad +e^{\ell t} \left(e^{2 t}-1\right)
(\ell^2-1)y_0 \Big), \vspace{0.2cm}\\
y(t)=& \dfrac{e^{-t}}{2 \left(\ell^2-1\right)} \Big(a
\left(-e^{(\ell+2) t} (\ell-1)^2+e^{\ell t} (\ell+1)^2-4 e^t
\ell\right)\\
& \qquad \qquad \quad +2 e^{\ell t+t} \left(\ell^2-1\right) y_0
(\cosh t-\ell \sinh t)\Big).
\end{array}
\end{equation}

The solution of system \eqref{+++} starting at the point
$(x,y)=(0,y_0)$ is given in \eqref{e2}.

\smallskip

Let $t_-$ and  $t_+$ be again the finite positive times defined in a
similar way to the Case 1, and if we have a periodic solution of
system \eqref{---}$+$\eqref{+++} one of the sets of three equations
\eqref{e2a}, or \eqref{e2b} must be satisfied. Now we shall assume
that equations \eqref{e2a} hold. Again the proof of Theorem \ref{t1}
in Case 2 if equations \eqref{e2b} hold is completely analogous to
the proof when equations \eqref{e2a} hold.

\smallskip

Using \eqref{e2}, \eqref{e20}, taking into account that $\ell\ne \pm
1$ and that the focus must be on the discontinuity straight line
(i.e. $c=0$) the three equations \eqref{e2a} become
\begin{equation*}\label{e3}
\begin{array}{rl}
e_1=& -e^{r t_+} (b+y_0) \sin t_+=0, \vspace{0.2cm}\\
e_2=& a \left(1-\ell-2 e^{\ell t_-+t_-}+e^{2 t_-}
(\ell+1)\right)+\left(e^{2 t_-}-1\right) (\ell^2-1)
y_0=0, \vspace{0.2cm}\\
e_3=& 2 (\ell^2-1) \left(-b+e^{r t_+} (b+y_0)(\cos t_+ - r \sin
t_+)\right)\\
& -e^{t_-} \Big(a e^{(-\ell-2) t_-} \left(-(\ell-1)^2+e^{2 t_-}
(\ell+1)^2-4 e^{\ell t_-+t_-} \ell\right)\\
& \qquad \quad +2 e^{(-\ell-1) t_-} (\ell^2-1) y_0 (\cosh t_-+\ell
\sinh t_-)\Big) =0.
\end{array}
\end{equation*}

{From} equation $e_1=0$, if $y_0 +b=0$, then at most there is one
periodic solution, and we are done. So we suppose that $\sin t_+=0$,
i.e. $t_+= \pi$. Now solving equation $e_2=0$ with respect to the
variable $y_0$ we get
\[
y_0= -\frac{a (1- \ell + e^{2t_-} + e^{2t_-} \ell- 2 e^{t_- + \ell
t_-})}{(e^{2t_-}-1) (\ell^2-1)}.
\]
Therefore the equation $e_3=0$ writes
\begin{equation*}\label{e21}
\begin{array}{l}
2 a e^{t_-} + \left(b \left(1+e^{\pi  r}\right)
\left(\ell^2-1\right)-a \left(e^{\pi  r}
(\ell-1)+\ell+1\right)\right)e^{\ell t_-}+\\
a \left(\ell+e^{\pi r} (\ell+1)-1\right)e^{(\ell+2) t_-}-\left(2
e^{\pi  r} a+b \left(1+e^{\pi r}\right) \left(\ell^2-1\right)\right)
e^{(2 \ell+1) t_-} =0,
\end{array}
\end{equation*}
or equivalently
\begin{equation*}\label{e22}
\begin{array}{l}
2 a f_0(t_-) + \left(b \left(1+e^{\pi  r}\right)
\left(\ell^2-1\right)-a \left(e^{\pi  r}
(\ell-1)+\ell+1\right)\right) f_1(t_-)+\\
a \left(\ell+e^{\pi r} (\ell+1)-1\right)f_2(t_-)-\left(2 e^{\pi  r}
a+b \left(1+e^{\pi r}\right) \left(\ell^2-1\right)\right) f_3(t_-)
=0,
\end{array}
\end{equation*}
where
\begin{equation*}\label{e23}
\begin{array}{l}
f_0(t_-)= e^{t_-}, \vspace{0.2cm}\\
f_1(t_-)= e^{\ell t_-}, \vspace{0.2cm}\\
f_2(t_-)= e^{(\ell+2) t_-}, \vspace{0.2cm}\\
f_3(t_-)= e^{(2 \ell+1) t_-}.
\end{array}
\end{equation*}

The set of functions $\{f_0,f_1,f_2,f_3\}$ is an ECT--system because
the Wronskians $W(f_0,...,f_k)(t_-)$ are not zero for $k=0,1,2,3$.
Indeed, we have
\[
\begin{array}{rl}
W(f_0)(t_-)=& e^{t_-} \neq 0, \vspace{0.2cm}\\
W(f_0,f_1)(t_-)=& (\ell-1) e^{(\ell+1) t_-}\neq 0, \vspace{0.2cm}\\
W(f_0,f_1,f_2)(t_-)=& 2(\ell^2-1)e^{(2 \ell+3) t_-} \neq 0, \vspace{0.2cm}\\
W(f_0,f_1,f_2,f_3)(t_-)=& 4\ell(\ell^2-1)^2e^{4 (\ell+1) t_-}\neq 0
\quad \mbox{if $\ell\neq 0$}.
\end{array}
\]
Therefore, if $\ell\neq 0$ the system \eqref{e2a} can have at most
$3$ solutions. Since the four coefficients of $f_0,f_1,f_2,f_3$ are
functionally dependent the system \eqref{---}$+$\eqref{+++} perhaps
do not reach the three solutions. But since the coefficients of
$f_1,f_2,f_3$ are functionally independent the system
\eqref{---}$+$\eqref{+++} can have two solutions. This completes the
proof of Theorem \ref{t1} in Case 3 when $\ell\neq 0$.

\smallskip

Assume $\ell=0$. Then the equation $e_3=0$ becomes
\begin{equation*}\label{e21}
\begin{array}{l}
a \left(e^{\pi  r}-1\right)-b \left(1+e^{\pi  r}\right)- \left(2a(
e^{\pi r}-1)-b \left(1+e^{\pi r}\right)\right)e^{t_-} \\
\qquad \qquad \qquad \qquad \qquad \qquad \qquad \qquad \,\,\,\,\,\,
+a \left(e^{\pi r}-1\right)e^{2 t_-} =0,
\end{array}
\end{equation*}
or equivalently
\begin{equation*}\label{e22}
\begin{array}{l}
(a \left(e^{\pi  r}-1\right)-b \left(1+e^{\pi  r}\right))f_0(t_-)-
\left(2a( e^{\pi r}-1)-b \left(1+e^{\pi r}\right)\right)f_1(t_-) \\
\qquad \qquad \qquad \qquad \qquad \qquad \qquad \qquad \qquad
\qquad \,\,\, +a \left(e^{\pi r}-1\right)f_2(t_-) =0,
\end{array}
\end{equation*}
where
\begin{equation*}\label{e23}
\begin{array}{l}
f_0(t_-)= 1, \vspace{0.2cm}\\
f_1(t_-)= e^{t_-}, \vspace{0.2cm}\\
f_2(t_-)= e^{2t_-}.
\end{array}
\end{equation*}

The set of functions $\{f_0,f_1,f_2\}$ is an ECT--system because the
Wronskians $W(f_0,...,f_k)(t_-)$ are not zero for $k=0,1,2$. Indeed,
we have
\[
\begin{array}{rl}
W(f_0)(t_-)=& 1 \neq 0, \vspace{0.2cm}\\
W(f_0,f_1)(t_-)=& e^{t_-}\neq 0, \vspace{0.2cm}\\
W(f_0,f_1,f_2)(t_-)=& 2e^{3t_-} \neq 0.
\end{array}
\]
Since the three coefficients of $f_0,f_1,f_2$ are functionally
dependent the system \eqref{---}$+$\eqref{+++} perhaps do not reach
the two solutions. But since the coefficients of $f_1,f_2$ are
functionally independent the system \eqref{---}$+$\eqref{+++} can
have one solution. This completes the proof of Theorem \ref{t1} in
Case 3 when $\ell= 0$.

\smallskip

In short, Theorem \ref{t1} is proved.

\section*{Acknowledgements}

The first author is supported by the FAPESP-BRAZIL grants 2010/
18015-6 and 2012/05635-1. The second author is partially supported
by a MCYT/FEDER grant MTM2008--03437, an AGAUR grant number
2014SGR--568, an ICREA Academia, two grants FP7-PEOPLE-2012-IRSES
316338 and 318999, and FEDER-UNAB10-4E-378.

\end{document}